\documentclass[reqno]{amsart}

\usepackage{amsmath}
\usepackage{amssymb}
\usepackage{enumerate}

\newtheorem{thm}{Theorem}
\newtheorem{lem}{Lemma}
\newtheorem{cor}{Corollary}
\newtheorem{propn}{Proposition}

\newtheorem{nonumthm}{Theorem}

\theoremstyle{definition}

\newtheorem{defn}{Definition}

\newtheorem{remark}{Remark}

\newcommand{\lemref}[1]{Lemma~\ref{#1}}
\newcommand{\thmref}[1]{Theorem~\ref{#1}}
\newcommand{\propref}[1]{Proposition~\ref{#1}}
\newcommand{\corref}[1]{Corollary~\ref{#1}}
\newcommand{\eqnref}[1]{(\ref{#1})}
\newcommand{\remref}[1]{Remark~\ref{#1}}
\newcommand{\agemo}{\mho}

\numberwithin{equation}{section}

\begin{document}

\title {The Coexponent of a Regular $p$-Group} 
\author{Paul J.\ Sanders}
\address
{Mathematics Institute\\ University of Warwick\\ Coventry\\ England CV4 7AL}
\email {pjs@maths.warwick.ac.uk}

\thanks
{Partially funded by EPSRC and the Mathematics Institute of the Universtity of
Warwick}

\baselineskip = 15pt

\begin{abstract}
A sharp bound is derived for the nilpotency class of a regular $p$-group
in terms of its coexponent, and is used to show that the number of groups of 
order $p^n$ with a given fixed coexponent, is independent of $n$, for $p$ and 
$n$ sufficiently large. Explicit formulae are calculated in the case of 
coexponent 3. 
\end{abstract}

\maketitle

\section{Introduction}
\begin{defn}
For a finite $p$-group $P$ of order $p^n$ and exponent $p^e$, we define its
coexponent $f(P)$ to be $n-e$.
\end{defn}

For a finite $p$-group $P$ of odd order, it is shown in \cite{boundpaper} that
${\rm cl}(P)\leq 2f(P)$ and that by excluding the prime $3$ as well, the
bound ${\rm cl}(P)\leq 2f(P)-1$ holds. The question of whether this bound is 
best possible was considered in \cite{coex}, but not resolved, and  
it was conjectured that for a non-negative integer $f$ it is
possible to exclude finitely many primes so that if $P$ is a finite
$p$-group of coexponent $f$, and $p$ is not one of the exluded primes, then 
$cl(P) \leq f+1$. In the first part of this paper we prove the following 
theorem which provides an affirmative solution to this conjecture.

\begin{thm}
\label{bound}
Let $P$ be a finite $p$-group where $p>f(P) + 1$. Then 
${\rm cl}(P) \leq f(P)+1$.
\end{thm}

\noindent Moreover, for each non-negative integer $f$ and prime $p\geq f+1$,
we construct an example of a finite $p$-group of coexponent $f$ whose
nilpotency class is exactly $f+1$.

In the second part of the paper we investigate the existence of stability for
the number of isomorphism classes of groups of order $p^n$ and coexponent $f$,
as $p$ and $n$ vary with $f$ an arbitrary fixed integer greater than $2$.
The works of Burnside \cite{burnside}, Miller \cite{miller} and others (see
\cite{miller} for a complete list), show that for coexponent $2$, the number of
isomorphism classes is constant for $p\geq 5$ and $n\geq 6$, and in this paper
we show the following result.

\begin{thm}
\label{stability}
For a fixed integer $f\geq 3$, the number $\Psi_{n,f}^p$ of isomorphism classes
of groups of order $p^n$ and coexponent $f$, depends only on $p$, for $p>f+1$ 
and $n\geq 3f$.
\end{thm}

The proof of this theorem is constructive and uses the bound on the nilpotency
class in \thmref{bound} to translate the problem into an equivalent question
about $p^n$-element Lie rings via the Magnus-Lazard Lie ring functors 
(\cite{magnus} and \cite{lazard}). In this setting, explicit calculations are
generally easier and in the final part of this paper we show how to derive
the following formulae for $\Psi_{n,3}^p$ where $p\geq 5$ and $n\geq 7$.

\begin{thm}
\label{formulae}
For $p\geq 5$ and $n\geq 7$ we have
\begin{equation*}
\Psi_{n,3}^p = 5p + 2\gcd(p-1,3) + \gcd(p-1,4) +
\left\{\begin{array}{ll}59,&n=7\\ 61,&n=8\\ 62,&n\geq 9\end{array}\right.
.
\end{equation*}
\end{thm}

This paper comprises an improved treatment of some of the material contained in
the author's doctoral dissertation \cite{mythesis}.

\section{Regular $p$-groups and the proof of\ \thmref{bound}}
For a finite $p$-group $P$ and each non-negative integer $i$, we define 
subgroups
$\Omega_i(P)=\langle x \in P : x^{p^i}=1 \rangle$ and 
$\agemo_i(P)=\langle x^{p^i} : x\in P\rangle$, and $P$ is called {\bf regular} 
if for any two elements $x,y\in P$ there exists 
$z\in \agemo_1(\gamma_2(\langle x,y \rangle))$ with $x^p y^p = (xy)^p z$. 
Denoting the exponent of a finite $p$-group $P$ by $p^{\mu(P)}$, we define
a sequence of non-negative integers 
$\omega_1(P),\omega_2(P),\ldots,\omega_{\mu(P)}(P)$ by the condition that
$p^{\omega_i(P)} = \vert \Omega_i(P) / \Omega_{i-1}(P) \vert$ for each
$i,\,1\leq i\leq \mu(P)$. For a regular $p$-group $P$ of order $p^n$, this 
sequence is non-increasing, forms a partition of $n$, and we denote the
dual partition by $\mu_1(P),\mu_2(P),\ldots,\mu_{\omega_1(P)}(P)$ and call this
sequence the {\it type invariants} of $P$. For a regular $p$-group $P$, there 
is a duality between the $\Omega$- and 
$\agemo$-series in that for each $i,\,1\leq i\leq \mu(P)$, we have
$\vert \Omega_i(P) / \Omega_{i-1}(P) \vert = \vert \agemo_{i-1}(P) / 
\agemo_i(P) \vert$. In particular, $\vert P / \agemo_1(P) \vert = 
p^{\omega_1(P)}$.

Regular $p$-groups were introduced by P.\ Hall in \cite{hall1} and
included a generalisation of the basis theorem for finite Abelian 
$p$-groups which we shall need. A sequence $(g_1,\ldots,g_r)$ of non-identity
elements of a
regular $p$-group $P$ is called a {\bf uniqueness basis} if each element of $P$
is uniquely expressible as $g_1^{k_1}\cdots g_r^{k_r}$ with 
$0\leq k_j < {\rm order}(g_j)$ for each $j,\,1\leq j\leq r$.

\begin{nonumthm}[P.\ Hall \mbox{\cite[\S4.5]{hall1}}]
Let $P$ be a finite regular $p$-group. Then 
\begin{enumerate}
\item If $(g_1,\ldots,g_r)$ is any uniqueness basis of $P$ then $r=\omega_1(P)$
and the set of orders of the elements of the basis is 
$\{p^{\mu_1(P)},\ldots,p^{\mu_{\omega_1(P)}(P)}\}$.
\item If $(g_1,\ldots,g_{\omega_1(P)})$ is any uniqueness basis of $P$ then any 
element $g_1^{k_1}\cdots g_{\omega_1(P)}^{k_{\omega_1(P)}}$ belongs to
$\agemo_i(P)$ if and only if each factor $g_j^{k_j}$ belongs to $\agemo_i(P)$.
\item $P$ possesses a uniqueness basis $(g_1,\ldots,g_{\omega_1(P)})$ and a
chain $P=N_0 \supseteq N_1 \supseteq \cdots \supseteq N_{\omega_1(P)} = 
\agemo_1(P)$ of normal subgroups with the property that for each
$i=1,\ldots,\omega_1(P)$, the element $g_i$ has order $p^{\mu_i(P)}$ and is
contained in  $N_{i-1}\setminus N_i$. 
\end{enumerate}
\end{nonumthm}

We now present the proof of \thmref{bound}. First observe that if $P$ is any 
finite $p$-group and $Q$ is a cyclic subgroup of index $p^{f(P)}$, then 
$Q^p \leq \agemo_1(P)$, and so $\vert P/\agemo_1(P) \vert \leq p^{f(P)+1}$. 
Therefore, if $p > f(P) + 1$, we have $\vert P/\agemo_1(P) \vert < p$, and 
P.\ Hall showed in \cite{hall2} that this condition implies that $P$ must be 
regular.
Hence to prove \thmref{bound} it suffices to show the following result.

\begin{propn}
\label{regbound}
Let $P$ be a finite regular $p$-group. Then ${\rm cl}(P) \leq f(P) + 1$.
\end{propn}

\begin{proof}
The case $f=0$ is clear, so we may assume that $f\geq 1$ i.e.\ 
$\omega_1(P) > 1$. Let $(g_1,\ldots,g_{\omega_1(P)})$ be a uniqueness basis of
$P$ where $g_i$ has order $p^{\mu_i(P)}$, and consider the
element $g^\prime=g_1^{p^{\mu_2(P)}}$. By regularity, each commutator 
$[g_1,g_i]$ has order at most $p^{\mu_2(P)}$ and so this implies that 
$g^\prime\in Z(P)$. Now the quotient 
${\overline P}=P/\langle g^\prime \rangle$ is a regular $p$-group of order 
$\vert P \vert / p^{\mu_1(P)-\mu_2(P)}$, and for $2\leq i\leq \omega_1(P)$,
the order of $g_i\langle g^\prime \rangle$ is $p^{\mu_i(P)}$ (because
$\langle g_i \rangle \cap \langle g_1 \rangle = 1$). It follows that
$(g_1\langle g^\prime \rangle,\ldots,g_{\omega_1(P)}\langle g^\prime \rangle)$
is a uniqueness basis of $\overline P$ and therefore the type invariants of 
$\overline P$ are $(\mu_2(P),\mu_2(P),\mu_3(P),\ldots,\mu_{\omega_1(P)}(P))$. 
Hence the exponent of $\overline P$ is $p^{\mu_2(P)}$ and so $\overline P$
has coexponent $f$. Therefore by passing to a central quotient we can assume
that $\mu_1(P) = \mu_2(P)$ and the result will follow if we can show that 
$cl(P) \leq f$ in this situation. 

So let $(g_1,g_2,\ldots,g_{\omega_1(P)})$ be a uniqueness basis of $P$ asserted
by part $3$ of Hall's theorem stated above. Then there exists a normal subgroup
$N$ of $P$ ($N=N_2$ in the above statement) with the following properties

\begin{enumerate}
\item $\vert P : N \vert = p^2$ and $\agemo_1(P) \leq N$.
\item $g_1,g_2 \not\in N$ and $g_i\in N$ for $3\leq i\leq \omega_3(P)$.
\end{enumerate}
 
\noindent Now defining the subset $K$ to be 
$\langle g_1^p \rangle \langle g_2^p \rangle \langle g_3 \rangle \cdots
\langle g_{\omega_1(P)} \rangle$, we see that $k\subseteq N$, and
since $(g_1,\ldots,g_{\omega_1(P)})$ is a uniqueness basis we have
$\vert K \vert = \vert P \vert / p^2$ so that $K=N$. 
Now if $\mu_1(P)=1$ it is easy to see that $cl(P)\leq f$, so
assume that $\mu_1(P) > 1$, i.e.\ $g_1^p,g_2^p \neq 1$. Then $N$ is a regular
$p$-group with a uniqueness basis $(g_1^p,g_2^p,g_3,\ldots,g_{\omega_1(P)})$
(note that the orders in this sequence are not necessarily non-increasing), and
we claim that for $1\leq i\leq \mu_1(P) - 1$ the following inclusion holds
\begin{equation}
\label{inclusionclaim}
[\agemo_i(P),\underbrace{P,\ldots,P}_{\omega_{i+1}(P) - 1\atop{\rm times}}] 
\subseteq \agemo_{i+1}(P)
\end{equation} 

\noindent To see this, fix $1\leq i\leq \mu_1(P) - 1$ and first consider the
following decompositions (using part $2$ of Hall's theorem)
\begin{eqnarray}
\agemo_i(P) &=& \langle g_1^{p^i} \rangle \langle g_2^{p^i} \rangle 
\langle g_3^{p^i} \rangle \cdots \langle g_{\omega_1(P)}^{p^i} \rangle
\nonumber\\
\agemo_i(N) &=& \langle g_1^{p^{i+1}} \rangle \langle g_2^{p^{i+1}} \rangle
\langle g_3^{p^i} \rangle \cdots \langle g_{\omega_1(P)}^{p^i} \rangle
\nonumber\\
\agemo_{i+1}(P) &=& \langle g_1^{p^{i+1}} \rangle \langle g_2^{p^{i+1}} \rangle
\langle g_3^{p^{i+1}} \rangle \cdots \langle g_{\omega_1(P)}^{p^{i+1}} \rangle
\nonumber
\end{eqnarray}

\noindent Therefore $\agemo_i(P) \geq \agemo_i(N) \geq \agemo_{i+1}(P)$, and 
since $g_1^{p^i},g_2^{p^i}\neq 1$ we have 
$\vert \agemo_i(P) : \agemo_i(N) \vert = p^2$ so that 
\begin{equation}
\label{index}
\big\vert \agemo_i(P) / \agemo_{i+1}(P) : \agemo_i(N) / \agemo_{i+1}(P) 
\big\vert = p^2
\end{equation}

\noindent Now for any elements $x,y\in P$, the commutator $[x,y]$ lies in $N$,
and so since $P/\agemo_i(N)$ is regular, the commutator $[x^{p^i},y]$ lies in
$\agemo_i(N)$. Therefore $[\agemo_i(P),P] \subseteq \agemo_i(N)$ so that
$[\agemo_i(P)/\agemo_{i+1}(P), P/\agemo_{i+1}(P)] \subseteq \agemo_i(N)/
\agemo_{i+1}(P)$, and by using \eqref{index} we have

$$\big\vert \agemo_i(P)/\agemo_{i+1}(P) : \big[\agemo_i(P)/\agemo_{i+1}(P), 
P/\agemo_{i+1}(P)\big] \big\vert \geq p^2.$$

\noindent Now since $\vert \agemo_i(P)/\agemo_{i+1}\vert = p^{\omega_{i+1}}$,
an easy induction shows that for any $j,1\leq j\leq \omega_{i+1}(P)-1$,  

$$\big\vert \agemo_i(P)/\agemo_{i+1}(P) : \big[\agemo_i(P)/\agemo_{i+1}(P),P
\underbrace{/\agemo_{i+1}(P),\ldots,P/}_{j\ {\rm times}}\agemo_{i+1}(P)\big] 
\big\vert \geq p^{j+1},$$

\noindent and then the claim \eqref{inclusionclaim} follows by taking 
$j=\omega_{i+1}(P)-1$. 

Now since $\vert P/\agemo_1(P) \vert = \omega_1(P) \geq 2$ we have 
$\gamma_{\omega_1(P)}(P) \subseteq \agemo_1(P)$, and then using 
\eqref{inclusionclaim} we see that for any $i,1\leq i\leq \mu_1(P) - 1$,

$$\gamma_{\omega_1(P) + (\omega_2(P) - 1) + \cdots + (\omega_{i+1}(P) - 1)}(P)
\subseteq \agemo_{i+1}(P).$$ 

\noindent Therefore, taking $i=\mu_1(P) - 1$ and writing, we see that 
$\vert P \vert = p^n$ this implies that $\gamma_{n-\mu_1(P)+1}(P) = 1$, and so
$cl(P) \leq n - \mu_1(P) = f(P)$ as required.

\end{proof}

We conclude this section by constructing for any prime $p$ and positive integer
$f$ with 
$p\geq f+1$, a finite $p$-group of coexponent $f$ and nilpotency class
exactly $f+1$.

So let $p\geq f+1\geq2$ and denote by $A$ an elementary Abelian $p$-group of
order $p^{f+1}$. Now choose a basis $(x_1,\ldots,x_{f+1})$ of $A$ and consider 
the automorphism $\alpha$ of $A$ which
fixes $x_{f+1}$ and maps $x_i$ to $x_i x_{i+1}$ for $1\leq i\leq f$. 
Since $p\geq f+1$, $\alpha$ has order $p$, and so letting $g$ be a generator of
a cyclic group of order $p^{n-f}$ where $n\geq f+2$, we can form the 
semi-direct product $P=A\rtimes \langle g \rangle$ where $g$ acts by $\alpha$. 
$P$ therefore has order $p^{n+1}$ and for $2\leq i\leq f+2$ it is not hard to
see that $\gamma_i(P) = \langle x_i,\ldots,x_{f+1} \rangle$.

\begin{lem}
Let $a\in A$ and $k\in {\mathbb Z}$. Then $(g^k a)^{p^2} = g^{kp^2}$.
\end{lem}
 
\begin{proof}
By Hall's commutator collecting process in \cite{hall1},\ 
$g^{kp}a^p = (g^k a)^p c_2^{e_2} \cdots c_p^{e_p}$ where $c_i\in \gamma_i(P)$ 
and $e_i$ is the $i$th binomial coefficient for $2\leq i\leq p$. Therefore 
$g^{kp}a^p = (g^k a)^p c_p$ (since $\gamma_2(P)$ has exponent $p$), and then 
the result follows from the fact that 
$c_p\in \gamma_p(P) \subseteq \langle x_{f+1} \rangle \subseteq Z(P)$.
\end{proof}

From this lemma we see that the exponent of $P$ equals the order of $g$ and
so $P$ has coexpoenent $f+1$. Now observe that since $g$ acts by an 
automorphism of order $p$ and $n-f\geq 2$, the element $g^{p^{n-f-1}}$ has 
order $p$ and lies in $Z(P)$. Therefore 
$\langle g^{p^{n-f-1}}x_{f+1}^{-1} \rangle \subseteq Z(P)$ and the quotient
$P/\langle g^{p^{n-f-1}}x_{f+1}^{-1} \rangle$ is a group of order $p^n$ which
has coexponent $f$ and nilpotency class $f+1$, as required. 

\section{A Lie Ring correspondence and the proof of \thmref{stability}}

As a consequence of \thmref{bound}, any finite $p$-group $P$ with $p>f(P) + 1$
has nilpotency class less than $p$. This condition enables the group to be 
endowed with a Lie ring stucture in which certain 
structural properties of the group coincide with the obvious structural 
analogue in the Lie ring. Questions involving properties preserved under this
passage to the Lie ring therefore have equivalent formulations in terms of
Lie rings and can sometimes be easier to handle in this setting. The statement 
of \thmref{stability} is an example of such a question and we prove this theorem
by solving the Lie ring formulation. 

The Lie ring correspondence which we use was first discovered by Magnus in 
\cite{magnus} and later, independently, by Lazard in \cite{lazard}. For a prime
$p$, let $\Gamma_p$ denote the category of finite $p$-groups whose nilpotency
class is less than $p$, and let $\Lambda_p$ denote the category of finite
nilpotent Lie rings whose order is a power of $p$ and whose nilpotency class
is less than $p$. The main properties of the correspondence are summarised in
the following theorem. 

\begin{nonumthm}[Magnus \cite{magnus}, Lazard \cite{lazard}] 
Let $p$ be a prime, $P$ a group in $\Gamma_p$ and $L$ a Lie ring in $\Lambda_p$.
Then there exists a Lie ring ${\mathcal L}_p(P)$ in $\Lambda_p$ and a $p$-group 
${\mathcal G}_p(L)$ in $\Gamma_p$ such that 
\begin{enumerate}
\item $L$ has the same underlying set as ${\mathcal G}_p(L)$.
\item $P$ has the same underlying set as ${\mathcal L}_p(P)$.
\item\label{inversea} The Lie ring ${\mathcal L}_p({\mathcal G}_p(L))$ 
coincides with $L$. 
\item\label{inverseb} The $p$-group ${\mathcal G}_p({\mathcal L}_p(P))$ 
coincides with $P$.
\item\label{covariant} ${\mathcal G}_p$ and ${\mathcal L}_p$ are covariant 
functors between the respective categories.
\item\label{orderscoincide} The order of an element $g\in P$ coincides with its
additive order in
${\mathcal L}_p(P)$.
\item Subgroups of $P$ coincide with Lie-subrings of ${\mathcal L}_p(P)$, and
normal subgroups of $P$ coincide with ideals of ${\mathcal L}_p(P)$.
\item For normal subgroups $H$ and $K$ of $P$, the commutator group $[H,K]_P$
coincides with the ideal of ${\mathcal L}_p(P)$ generated by all Lie brackets
of elements of $H$ with elements of $K$. 
\end{enumerate}
\end{nonumthm}

\noindent
Part~\ref{orderscoincide} of this theorem implies that the exponent of $P$ 
equals
the exponent of the underlying Abelian group of ${\mathcal L}_p(P)$, and so 
the coexponents of these two $p$-groups also coincide. Parts~\ref{inversea},
\ref{inverseb}~and~\ref{covariant} of this theorem
imply that for a fixed natural number $n$ and prime $p$, the number of 
isomorphism classes of groups of order $p^n$ in $\Gamma_p$ equals the number
of isomorphism classes of $p^n$-element Lie rings in $\Lambda_p$. Therefore
\thmref{stability} is equivalent to the following theorem.

\begin{thm}
\label{formulation}
For a fixed integer $f\geq 3$, the number of isomorphism classes of 
$p^n$-element Lie rings in $\Lambda_p$ with additive coexponent $f$, depends
only on $p$, for $p > f + 1$ and $n \geq 3f$.
\end{thm} 

\noindent For the remainder of this section we concentrate on the proof of this
theorem. 

We first recall some definitions.
A derivation of a Lie ring $L$ is an Abelian group endomorphism 
$\phi : L \rightarrow L$ where $[x,y]\phi = [x\phi,y] + [x,y\phi]$ for each
$x,y\in L$. The derivations ${\rm Der}(L)$ form a Lie ring with bracket given
by $[\phi_1,\phi_2] = \phi_1\phi_2 -\phi_2\phi_1$. Each element $x\in L$ 
defines a derivation ${\rm ad}_{\rm r}x$ where 
${\rm ad}_{\rm r}x : y \mapsto [y,x]$. 
A derivation $\phi$ is said to be nilpotent if $\phi^k$ is identically zero
for some $k$, and is said to centralize $x\in L$ if $x\phi = 0$. We denote by
${\rm Der}(L)_x$ the Lie subring of ${\rm Der}(L)$ consisting of those
derivations centralizing $x$.
Given Lie rings $L,M$ and a Lie ring homomorphism 
$\theta : M \rightarrow {\rm Der}(L)$, we denote by $L \rtimes_\theta M$ the Lie
ring on $L \oplus M$ with bracket given by 
$[l_1 + m_1,l_2 + m_2] = [l_1,l_2]_L + [m_1,m_2]_M + l_1(m_2\theta) -
l_2(m_1\theta)$. $L \rtimes_\theta M$ satisfies an obvious universal property.

Now fix an arbitrary prime number $p$ and an arbitrary positive integer $f$, and
let $L$ be a nilpotent $p^n$-element Lie ring of additive coexponent $f$. Then,
as an Abelian group, $L=\langle u \rangle \oplus V$ where $u$ has order 
$p^{n-f}$ and $\vert V \vert =
p^f$. Observe that $\mu_2(L) = \mu_1(V)\leq f$, and if 
$n - f \geq 2\mu_2(L)$ then 
$p^{n-f - \mu_2(L)}u \in Z(L)\cap\Omega_{\mu_2(L)}(L)$. In this situation,
$\Omega_{\mu_2(L)}(L)=\langle p^{n-f - \mu_2(L)}u \rangle \oplus V$ is a 
nilpotent ideal of order $p^{\mu_2(L) + f}$ and the
derivation $\phi = {\rm ad}_r u |_{\Omega_{\mu_2(L)}(L)}$ is nilpotent and
centralizes a central element of order $p^{\mu_2(L)}$ (namely 
$p^{n-f - \mu_2(L)}u$).
It is then easy to see that for any generator $w$ of the Abelian Lie
ring ${\mathbb Z}/p^{n-f}{\mathbb Z}$, there is a Lie ring epimorphism 
$\psi : \Omega_{\mu_2(L)} \rtimes_\theta ({\mathbb Z}/p^{n-f}{\mathbb Z})w
\twoheadrightarrow L$ where $w\theta = \phi$ and 
${\rm Ker}(\psi) = \langle p^{n-f - \mu_2(L)}u - w \rangle$.
Observe that $\Omega_{\mu_2(L)}(L)$ is a Lie ring $U$ satisfying the following 
properties

\begin{enumerate}[{P}-1]
\item $U$ is nilpotent, has additive coexponent $f$ and 
$\mu_1(U) = \mu_2(U) \leq f$.\label{Uone} 
\item \label{Utwo} $\vert U \vert = p^{f+\mu_1(U)}$. 
\item \label{Uthree} $U$ possesses a central element of order 
$p^{\mu_1(U)}$.
\end{enumerate}

\begin{lem}
\label{construct}
Let $U$ satisfy {\rm P}-$\ref{Uone}$ -- {\rm P}-$\ref{Uthree}$ and $z$ be any 
element of $U$ asserted by {\rm P}-$\ref{Uthree}$. Then for any Lie ring 
homomorphism  
$\theta : ({\mathbb Z}/p^{n-f}{\mathbb Z})w \rightarrow {\rm Der}(U)_z$ 
where $w$ is a generator of ${\mathbb Z}/p^{n-f}{\mathbb Z}$ and 
$n-f \geq 2\mu_1(U)$, we have 
\begin{enumerate}
\item The subgroup $\langle p^{n-f-\mu_1(U)} w - z \rangle$ is an ideal of
$U \rtimes_\theta ({\mathbb Z}/p^{n-f}{\mathbb Z})w$ \label{ideal}.
\item The quotient $U(n-f,w,w\theta,z) = U \rtimes_\theta 
({\mathbb Z}/p^{n-f}{\mathbb Z})w / \langle p^{n-f-\mu_1(U)} w - z \rangle$ is 
a Lie ring of order $p^n$
and additive coexponent $f$.\label{quotient}
\item $\mu_2(U(n-f,w,w\theta,z)=\mu_1(U)$ and 
$\Omega_{\mu_1(U)}(U(n-f,w,w\theta,z)) \cong U$\label{Uembeds}.
\item If $w\theta$ is nilpotent then so is $U(n-f,w,w\theta,z)$
\label{nilpotent}.
\end{enumerate} 
\end{lem}

\begin{proof}
The condition $n-f- \mu_1(U)\geq \mu_1(U)$ ensures that 
$p^{n-f- \mu_1(U)}w\in {\rm Ker}(\theta)$, and Part \ref{ideal} follows easily 
from this, $z\in Z(U)$ and $w\theta$ centralizing $z$. For part \ref{quotient},
the order of $\langle p^{n-f-\mu_1(U)} w - z \rangle$ is $p^{\mu_1(U)}$ and so 
$U(n-f,w,w\theta,z)$ is a Lie 
ring of order $p^n$ by property P-\ref{Utwo}. Moreover,
$\langle w \rangle \cap \langle p^{n-f-\mu_1(U)} w - z \rangle = \{0\}$, and
so $U(n-f,w,w\theta,z)$ has coexponent $f$. 
For Part \ref{Uembeds}, note that 
$U \cap \langle p^{n-f-\mu_1(U)} w - z \rangle = \{0\}$
and then it is straightforward to see that $U$ maps isomorphically onto
$\Omega_{\mu_2(U)}(U(n-f,w,w\theta,z))$. 
To show part \ref{nilpotent}, we denote the $k$-th term of the lower central
series of a Lie ring $M$ by $M_k$ and claim that for any $j\geq 1$ there
is a natural number $n_j$ such that $U(n-f,w,w\theta,z)_{n_j} \subseteq U_j$. 
The case $j=1$ is
clear since $U(n-f,w,w\theta,z)/U$ is Abelian, so suppose we have shown it for 
some $j\geq 1$.
Then 
$U(n-f,w,w\theta,z)_{n_j+1} 
\subseteq I_{j+1}+U(n-f,w,w\theta,z)_{n_j}{\rm ad_r}w$.
Using this as the base
case for a second induction, it is straightforward to show that
$U(n-f,w,w\theta,z)_{n_j + k} \subseteq U_{j+1} + 
U(n-f,w,w\theta,z)_{n_j}({\rm ad_r} w)^k$ for any $k\geq 1$.
Nilpotency of ${\rm ad_r}w$ then implies that 
$U(n-f,w,w\theta,z)_{n_j+k}\subseteq U_{j+1}$
for some $k\in {\mathbb N}$, and so $U(n-f,w,w\theta,z)$ is nilpotent. 
\end{proof} 

Continuing to use the notation of the previous lemma, we have the following
criteria for isomorphism 

\begin{lem}
\label{criteria}
$U(n-f,w,w\theta,z) \cong U(n-f,w^\prime,w^\prime\theta^\prime,z^\prime)$
if and only if there exists a Lie ring automorphism $\pi$ of $U$ with
\begin{enumerate}
\item $z \pi = \alpha z^\prime$ for some $\alpha\in {\mathbb Z}$ coprime to $p$.
\item $\pi^{-1} (w\theta) \pi = {\rm ad}_rx + \alpha(w^\prime\theta^\prime)$ for
some $x\in U$. 
\end{enumerate}
\end{lem} 

\begin{proof}
Suppose first that there is a Lie ring isomorphism 
$\eta : U(n-f,w,w\theta,z) 
\rightarrow U(n-f,w^\prime,w^\prime\theta^\prime,z^\prime)$. By part 
\ref{Uembeds} 
of \lemref{construct}, $U$ can be naturally identified with the ideals 
$\Omega_{\mu_1(U)}(U(n-f,w,w\theta,z))$ and 
$\Omega_{\mu_1(U)}(U(n-f,w^\prime,w^\prime\theta^\prime,z^\prime))$, and 
$w,w^\prime$ both have order 
$p^{n-f}$ in $U(n-f,w,w\theta,z),
U(n-f,w^\prime,w^\prime\theta^\prime,z^\prime)$ 
respectively. Under these identifications, $U$ is invariant under 
$\eta$ so that $\eta \vert_U$ is a Lie ring automorphism, and for $w\eta$ to
have order $p^{n-f}$ we must have $w\eta = x + \alpha w^\prime$ where $x\in U$ 
and $(\alpha,p) = 1$. Therefore $z\eta = (p^{n-f-\mu_1(U)}w) \eta =
p^{n-f-\mu_1(U)}(x + \alpha w^\prime) = \alpha z^\prime$ 
(since $n-f\geq2\mu_1(U)$). Now denoting the Lie ring brackets on
$U(n-f,w,w\theta,z)$ and $U(n-f,w^\prime,w^\prime\theta^\prime,z^\prime)$ by
$[.,.]$ and $[.,.]^\prime$ respectively, it follows that for any $u\in U$ we 
have $u(\eta|_U^{-1} (w\theta^\prime) \eta|_U) =
[u\eta|_U^{-1},w]\eta = [u,w\eta]^\prime = u({\rm ad_r}x +
\alpha (w^\prime\theta^\prime))$, and so conditions 1 and 2 hold with 
$\pi = \eta|_U$.
Conversely, suppose that conditions 1 and 2 hold, and define a homomorphism
of Abelian groups 
$\chi : U \oplus ({\mathbb Z}/p^{n-f}{\mathbb Z})w \rightarrow
U \oplus ({\mathbb Z}/p^{n-f}{\mathbb Z})w^\prime$ where $\chi|_U = \pi$ and
$w\chi = x + \alpha w^\prime$. Since $(\alpha,p) = 1$, $\chi$ is bijective, and
from condition 2 it is straightforward to see that $\chi$ is a Lie ring
isomorphism $U \rtimes_{\theta} ({\mathbb Z}/p^{n-f}{\mathbb Z})w
\rightarrow U \rtimes_{\theta^\prime} 
({\mathbb Z}/p^{n-f}{\mathbb Z})w^\prime$. 
Now letting $I = \langle z - p^{n-f-\mu_1}w \rangle$ and
$I^\prime = \langle z^\prime - p^{n-f-\mu_1}w^\prime \rangle$, it is
straightforward to see that these are both central ideals of their
respective Lie rings, and $I \chi = I^\prime$. Therefore $\chi$ induces the 
required isomorphism.
\end{proof}

\begin{propn}
\label{independent}
Let $(\lambda_1,\ldots,\lambda_k)$ be a partition of $f$ with
$\lambda_1 \geq \cdots \geq \lambda_k$, and suppose that $n \geq f+2\lambda_1$.
Then for any positive integer $c$, the number of isomorphism classes of
nilpotent $p^n$-element Lie rings $L$ where $L$ has nilpotency class $c$,
$\mu_2(L)=\lambda_1$ and $\Omega_{\mu_2(L)}(L)$ has type invariants
$(\lambda_1,\lambda_1,\ldots,\lambda_k)$, is independent of $n$.
\end{propn}

\begin{proof}
Suppose that $U$ is a Lie ring satisfying 
{\rm P}-$\ref{Uone}$ -- {\rm P}-$\ref{Uthree}$ and whose type
invariants are $(\lambda_1,\lambda_1,\ldots,\lambda_k)$. 
Observe that the additive exponent of ${\rm Der}(U)$ is a divisor of
the additive exponent of $U$, namely $p^{\lambda_1}$. Therefore, choosing 
$w_{n-f}\in({\mathbb Z}/p^{n-f}{\mathbb Z})^\times$ arbitrarily, we see that 
for any $\sigma\in{\rm Der}(U)$ there is a Lie ring homomorphism
$\theta : ({\mathbb Z}/p^{n-f}{\mathbb Z})w_{n-f} \rightarrow {\rm Der}(U)$ with
$w\theta = \sigma$. Then from
\lemref{construct}, the discussion preceding it, and \lemref{criteria}, we
see that for $n\geq f+2\lambda_1$, 
the number of isomorphism classes of nilpotent $p^n$-element Lie rings $L$
of additive coexponent $f$ with $\mu_2(L)=\lambda_1$ and 
$\Omega_{\mu_2(L)}\cong U$, equals the number of isomorphism classes arising 
from the (finite) set of Lie rings
\begin{eqnarray}
\label{reps}
\lefteqn{\mathcal{N}(U,n-f,w_{n-f})=}\nonumber\\
&\{U(n-f,w_{n-f},\sigma,z) : z\in Z(U)\ {\rm of\ order}\ p^{\lambda_1},
 \sigma\in{\rm Der}(U)_z\ {\rm is\ nilpotent} \}
\end{eqnarray}

\noindent Moreover, it is straightforward to see that for 
$n,n^\prime\geq f+2\lambda_1$, the bijection
\begin{eqnarray*}
\mathcal{N}(U,n-f,w_{n-f}) &\rightarrow& 
\mathcal{N}(U,n^\prime-f,w_{n^\prime-f})\\
U(n-f,w_{n-f},\sigma,z) &\mapsto& U(n^\prime-f,w_{n^\prime-f},\sigma,z)
\end{eqnarray*}
preserves the nilpotency class. Letting $U$ vary over a transversal for the
isomorphism classes of Lie rings satisfying {\rm P}-$\ref{Uone}$ -- 
{\rm P}-$\ref{Uthree}$ and which have type invariants 
$(\lambda_1,\lambda_1,\ldots,\lambda_k)$, gives the required result. 

\end{proof}

\noindent We can now prove the main theorem of this section.

\begin{proof}[Proof of \thmref{formulation}]
For $n\geq 3f$, proposition \ref{independent} implies that 
given a partition $(\lambda_1,\ldots,\lambda_k)$ of $f$ with 
$\lambda_1 \geq \cdots \geq \lambda_k$, the number of nilpotent $p^n$-element
Lie rings $L$ where $\mu_2(L)=\lambda_1$, $\Omega_{\mu_2(L)}(L)$ has type 
invariants $(\lambda_1,\lambda_1,\ldots,\lambda_k)$, and $L$ has nilpotency
class less than $p$, is independent of $n$. The theorem follows by taking this 
result over all partitions of $f$. 
\end{proof}

\section{Formulae for the number of groups of order $p^n$ and coexponent 3 
($p\geq 5,n\geq 7$)}

In this section we show how to use the Lie ring results of the previous section
to derive the expressions given in the statement of \thmref{formulae} for 
$\Psi_{n,3}^p$ when $p\geq5$ and $n\geq 7$. We first need some notation. For a 
prime $p$ and positive integers $e \geq \lambda_1 \geq \cdots \geq \lambda_k$, 
we denote by $\Psi_e^p(\lambda_1,\ldots,\lambda_k)$, the number of isomorphism 
classes of regular $p$-groups of type $(e,\lambda_1,\ldots,\lambda_k)$. If
we now let $f$ be a positive integer and $p$ a prime with $p>f+1$, then as a 
consequence of \propref{regbound}, any $p$-group of coexponent $f$ is regular, 
and so for $n\geq 2f$, we have 

\begin{equation}
\label{partitionsum}
\Psi_{n,f}^p = \sum_{{\rm partitions} \atop{\underline{\lambda}\ {\rm of\ f}}} 
\Psi_{n-f}^p(\underline{\lambda})
\end{equation}

\noindent The formulae given in \thmref{formulae} follow immediately from 
\eqnref{partitionsum} and the following proposition.
\begin{propn}
\label{decomposition}
Let $p$ be a prime greater than or equal to 5. Then
\begin{enumerate}
\item $\Psi_{n-3}^p(1,1,1) = 23 + 2\gcd(p-1,3) + \gcd(p-1,4), \quad n\geq 5$.
\label{usep5}
\item $\Psi_{n-3}^p(2,1) = 5p + 30, \quad n\geq 7.$ 
\label{laborious}
\item $\Psi_{n-3}^p(3) = \left\{\begin{array}{lr}
6,&\mbox{$n=7$}\\ 8,&\mbox{$n=8$}\\ 9,&\mbox{$n\geq 9$}\end{array}\right.$
\end{enumerate}
\end{propn}

A determination of the $p$-groups of coexponent 3 has been previously attempted
with a solution for $p=2$ appearing in \cite{mckelden}, and a solution for
odd primes $p$ being claimed in \cite{neikirk}. For $p\geq 5$, our formula
for $\Psi_{n-3}^p(3)$ agrees with that given in \cite{neikirk}, but our 
formulae for $\Psi_{n-3}^p(2,1)$ and $\Psi_{n-3}^p(1,1,1)$ do not agree  
(in \cite{neikirk}, the claimed formulae are $\Psi_{n-3}^p(2,1) = 5p + 32$
and $\Psi_{n-3}^p(1,1,1) = 23$ for $n\geq 7$). The approach used in 
\cite{neikirk} is via generators and relations, but there is no check that the 
groups presented have the correct order. For example, the first group in the
summary
table for $\Psi_{n-3}^p(2,1)$ is given by the presentation
\begin{equation*}
G = \langle u,v,w\ \vert\ u^{p^{n-3}} = v^{p^2} = w^p = 1,\ [u,v] = u^{p^{n-5}},
\ [u,w] = u^{p^{n-4}},\ [v,w] = v^p \rangle
\end{equation*}
and in this group we see that $u$ generates a normal cyclic subgroup $N$ which
must have order at most $p^{n-4}$ so that the automorphism induced on $N$ by
$v$ has order dividing $p$ (this is the case since $[v,w] = v^p$ and 
${\rm Aut}(N)$ is Abelian). But then $G$ is a split extension of the
cyclic group of order $p^{n-4}$ by the non-Abelian group of order $p^3$ and
exponent $p^2$. But then $G$ has order $p^{n-1}$ and not $p^n$ as claimed.

The only other paper we are aware of which considers specifically $p$-groups
of coexponent 3 is \cite{titov} in which the quotients arising modulo the core 
of a largest cyclic subgroup are determined.

We now show how the results of the previous section are used to derive the
formulae given in \propref{decomposition}. The following simple application
of the Magnus-Lazard Lie ring functors translates calculation of
$\Psi_e^p(\lambda_1,\ldots,\lambda_k)$ into a Lie ring setting.

\begin{propn}
\label{translate}
Let $p$ be a prime and $e\geq\lambda_1\geq \cdots \geq\lambda_k$ be positive
integers. Then for $p>\lambda_1 + \cdots \lambda_k + 1$, 
$\Psi_e^p(\lambda_1,\ldots,\lambda_k)$ equals the number of
isomorphism classes of nilpotent Lie rings of order 
$p^{e+\lambda_1+\cdots+\lambda_k}$ which have type invariants
$(e,\lambda_1,\ldots,\lambda_k)$ and nilpotency class less than $p$.
\end{propn}

\begin{proof}
From \propref{regbound}, the condition on $p$ means that any regular $p$-group
of type $(e,\lambda_1,\ldots,\lambda_k)$ belongs to $\Gamma_p$. The result
then follows by noting that since the functors ${\mathcal G}_p$ and 
${\mathcal L}_p$ preserve the order of an element, they must also preserve the
type invariants.
\end{proof}

\noindent We then have the following stability result.

\begin{cor}
\label{typestab}
Let $f$ be a positive integer and $p$ a prime with $p>f+1$. Then for any
partition $(\lambda_1,\ldots,\lambda_k)$ of $f$, 
$\Psi_{n-f}^p(\lambda_1,\ldots,\lambda_k)$ is independent of 
$n$ provided $n\geq f + 2\lambda_1$.
\end{cor}

\begin{proof}
This follows immediately from \propref{independent} and \propref{translate}.
\end{proof}

\subsection{Calculating $\Psi_{n-3}^p(1,1,1)$ for $p,n\geq5$.}\mbox{}

From \corref{typestab}, we see that $\Psi_{n-3}^p(1,1,1) = \Psi_2^p(1,1,1)$ 
for $p,n\geq5$, i.e.\ $\Psi_{n-3}^p(1,1,1)$ equals the number of groups of order
$p^5$ which have type invariants $(2,1,1,1)$. The most recently published
list of groups of order $p^5$ (that we know of) is contained in \cite{james}, 
and for 
$p\geq5$ these are tabulated by their type invariants. From this list we
obtain the formula given in part \ref{usep5} of \propref{decomposition}. Our
formula therefore depends on the correctness of the classification of groups of
order $p^5$ for $p\geq5$. These groups have been the subject of a number of
papers, among them \cite{bender}, \cite{schreier} and \cite{james}, and they
all agree on the correct formula for $p\geq5$. This formula is reiterated as
being correct in \cite{newman}. The dependence on the residue class
modulo 12 (not present in the formula for $\Psi_{n-3}^p(1,1,1)$ given in the 
paper \cite{neikirk}), is a contribution from the
groups of order $p^5$ and maximal nilpotency class 4, and this dependence is in
accordance with Blackburn's results in \cite{blackburn}. 

\subsection{Calculating $\Psi_{n-3}^p(2,1)$ for $p\geq5$ and $n\geq7$.}
\mbox{}
 
From \propref{translate} and \corref{typestab}, we see that for $p\geq5$ and 
$n\geq7$, $\Psi_{n-3}^p(2,1)$ equals the number of 
nilpotent Lie rings of order $p^7$ which have type invariants $(4,2,1)$ and
nilpotency class less than $p$. The proof of \propref{independent} 
gives the following procedure for calculating this number. 

\begin{enumerate}
\item
Find a transversal ${\mathcal U}_p$ for the nilpotent Lie rings of order
$p^5$ which have type $(2,2,1)$ and possess a central element of order $p^2$.
\item
For each Lie ring $U$ in ${\mathcal U}_p$, define an equivalence relation
$\sim$ on the set 
$\{(z,\sigma) : z\in Z(U)\ {\rm of\ order}\ p^2, \sigma\in {\rm Der}(U)_z
\ {\rm is\ nilpotent}\}$ by
$(z,\sigma) \sim (z^\prime,\sigma^\prime)$ if and only if there exists
$\pi\in {\rm Aut_{Lie}}(U)$ with 
$z\pi = \alpha z^\prime$ and $\sigma\pi = \pi({\rm ad_r}x + 
\alpha\sigma^\prime)$ for some $x\in U$ and $\alpha$ coprime to $p$.
\item
For each Lie ring $U$, calculate a complete set of representatives for the 
$\sim$-classes.
\item
For each representative $(z,\sigma)$, calculate the nilpotency class of
$U(4,w_4,\sigma,z)$. 
\item
The required formula is then obtained by evaluating, for
each prime $p$, the number of representatives which have nilpotency class less
than $p$. 
\end{enumerate} 

\begin{remark}
\label{lieringbound}
Step 4 is, in fact, unnecessary because the result analogous to 
\propref{regbound} holds for Lie rings. Namely, given a nilpotent Lie ring of 
order $p^n$ and additive coexponent $f$, the nilpotency class is at most
$f+1$. The proof of this is directly analogous to the proof of
\propref{regbound} and is given in detail in \cite{mythesis}. It follows that
since
we are only interested in $p\geq5$, the required formula can be calculated
after Step 3.
\end{remark}

To find a transversal ${\mathcal U}_p$, choose 
$z,u_1\in({\mathbb Z}/p^2{\mathbb Z})^\times$ and 
$u_2\in({\mathbb Z}/p{\mathbb Z})^\times$, and consider the Abelian group
$A=({\mathbb Z}/p^2{\mathbb Z})z \oplus ({\mathbb Z}/p^2{\mathbb Z})u_1\oplus
({\mathbb Z}/p{\mathbb Z})u_2$. It is straightforward to see that for any
element $v\in\Omega_1(A)$, there is a unique Lie ring structure on $A$ with
$z$ central and $[u_1,u_2] = v$, and nilpotency of such a Lie ring is
equivalent to $v$ having no non-zero component in $\langle u_2 \rangle$. 
Therefore, the nilpotent Lie ring structures on $A$ with $z$ central are given 
by $[u_1,u_2]=\alpha_1 p z + \alpha_2 p u_1$ where 
$\alpha_1,\alpha_2\in{\mathbb Z}$. Letting $V,W,X$ be the Lie rings where
$(\alpha_1,\alpha_2)= (0,1),(1,0),(0,0)$ respectively, we see that any Lie
ring on $A$ with $\alpha_2 \not\equiv0\pmod{p}$ is isomorphic to $V$, and any
Lie ring on $A$ with $\alpha_2\equiv0\pmod{p}$ is isomorphic to $W$ or $X$.
We can therefore take ${\mathcal U}_p = \{V,W,X\}$ ($V \not\cong W$ since
$[W,W]\cap\agemo_1(Z(W)) \neq \{0\}$ whereas 
$[V,V]\cap\agemo_1(Z(V)) = \{0\}$). 

In each of $V,W$ and $X$, the automorphism group acts transitively on the set
of central elements of order $p^2$. Therefore, in Step 2 of the above
procedure, we need only consider pairs $(z,\sigma)$ corresponding to a fixed
central element $z$ of order $p^2$, and so $\sim$ is an equivalence relation
on ${\rm Der}(U)_z$. We now consider Step 3 for each of $V,W,X$, and show
that the number of $\sim$-classes is $2p+1,3p+11,18$ respectively. By
\remref{lieringbound}, this then shows that the required formula holds.

\subsubsection{The Lie ring $V$.}
Let the basis of $V$ be ordered as $(z,u_1,u_2)$ so that we can
represent ${\rm Hom}_{\mathbb Z}(A,A)$ by the ring of matrices
\begin{equation}
\label{hom}
\left( 
\begin{array}{ccc}
{\mathbb Z}/p^2{\mathbb Z}&{\mathbb Z}/p^2{\mathbb Z}&{\mathbb Z}/p{\mathbb Z}\\
{\mathbb Z}/p^2{\mathbb Z}&{\mathbb Z}/p^2{\mathbb Z}&{\mathbb Z}/p{\mathbb Z}\\
p{\mathbb Z}/p^2{\mathbb Z}&p{\mathbb Z}/p^2{\mathbb Z}&{\mathbb Z}/p{\mathbb Z}
\end{array}
\right)
\end{equation}
Observe first that the nilpotent elements of ${\rm Der}(V)_z$ correspond to
certain matrices of the form
\begin{equation}
\label{nilpotentform}
\left(
\begin{array}{ccc}
0&0&0\\a_1&pa_2&a_3\\pb_1&pb_2&0
\end{array}
\right).
\end{equation}

\noindent In order that such a matrix $\kappa$ defines a derivation we need
only check that $[u_1,u_2]\kappa = [u_1\kappa,u_2] + [u_1,u_2\kappa]$, and
since $[u_1,u_2] = p u_1$ this condition holds if and only if 
$a_1 \equiv 0 \pmod{p}$. 
Now the inner derivations of $V$ are represented by the matrices 
\begin{equation*}
\left(
\begin{array}{ccc}
0&0&0\\
0&p{\mathbb Z}/p^2{\mathbb Z}&0\\
0&p{\mathbb Z}/p^2{\mathbb Z}&0
\end{array}
\right),
\end{equation*}
\noindent and so since $\kappa \sim \kappa + {\rm ad_r} x$ for any $x\in V$, it
follows that we need only calculate the number of $\sim$-classes among
nilpotent elements of ${\rm Der}(V)_z$ defined by the matrices 
\begin{equation}
\label{derV}
\left( 
\begin{array}{ccc}
0&0&0\\
p{\mathbb Z}/p^2{\mathbb Z}&0&{\mathbb Z}/p{\mathbb Z}\\
p{\mathbb Z}/p^2{\mathbb Z}&0&0
\end{array}
\right)
\end{equation}
\noindent Now let $\pi\in {\rm Hom}_{\mathbb Z}(V,V)$ fix $\langle z \rangle$ 
set-wise so that $\pi$ is represented as
\begin{equation}
\label{autform}
\left(
\begin{array}{ccc}
\alpha&0&0\\
\beta_1&\beta_2&\beta_3\\
p\gamma_1&p\gamma_2&\gamma_3
\end{array}
\right),
\end{equation}
\noindent and observe that $\pi\in {\rm Aut}_{\mathbb Z}(V)$ if and only if
$\alpha,\beta_2,\gamma_3$ are coprime to $p$. Since 
$\langle z \rangle \pi = \langle z \rangle$ and $z\in Z(V)$, the condition that
$\pi$ preserve
the Lie bracket is given by $[u_1,u_2] \pi = [u_1 \pi, u_2 \pi]$. This is
holds if and only if $p\beta_1 z+ p\beta_2 u_1 = p\beta_2\gamma_3 u_1$, i.e.\
if and only if $\beta_1 \equiv 0 \pmod{p}$ and $\gamma_3 \equiv 1 \pmod{p}$.
Therefore automorphisms of $V$ fixing $\langle z \rangle$ set-wise are of the
form
\begin{equation}
\left(
\begin{array}{ccc}
\alpha&0&0\\
p\beta_1&\beta_2&\beta_3\\
p\gamma_1&p\gamma_2&1
\end{array}
\right)
\quad{\rm with}\ \gcd(\alpha,p)=\gcd(\beta_2,p)=1.
\end{equation}
\noindent We now determine a system of congruences for $\sim$-equivalence
between two derivation $\sigma,\tau$ of the form \eqref{derV}. So let
\begin{equation*}
\sigma=\left(
\begin{array}{ccc} 
0&0&0\\
pa_1&0&a_3\\
pb_1&0&0
\end{array}
\right)
\quad{\rm and}\quad
\tau=\left(
\begin{array}{ccc}
0&0&0\\
pc_1&0&c_3\\
pd_1&0&0
\end{array}
\right).
\end{equation*}
\noindent Then $\sigma \sim \tau$ if and only if there exist matrices
\begin{equation*}
m = \left(
\begin{array}{ccc}
0&0&0\\
0&p\lambda_1&0\\
0&p\lambda_2&0
\end{array}
\right)
\quad{\rm and}\quad
\pi = \left(
\begin{array}{ccc}
\alpha&0&0\\
p\beta_1&\beta_2&\beta_3\\
p\gamma_1&p\gamma_2&1
\end{array}
\right)
\end{equation*}
\noindent with $\gcd(\alpha,p) = \gcd(\beta_2,p) = 1$, such that
$\sigma \pi = \pi (m + \alpha \tau)$. This is equivalent to the following
system of congruences having a solution for 
$\lambda_1,\lambda_2,\beta_1,\beta_2,\beta_3,\gamma_1,\gamma_2$
\begin{eqnarray*}
\alpha,\beta_2&\not\equiv&0\pmod{p}\\
p(a_1\alpha + a_3\gamma_1)&\equiv&p\alpha(\beta_2 c_1 + \beta_3 d_1)\pmod{p^2}
\\
p a_3 \gamma_2&\equiv&p(\beta_2\lambda_1 + \beta_3\lambda_2)\pmod{p^2}\\
a_3&\equiv&\beta_2\alpha c_3\pmod{p}\\
p b_1 \alpha&\equiv&p \alpha d_1\pmod{p^2}\\
p\lambda_2&\equiv&0\pmod{p^2}
\end{eqnarray*}
\noindent which reduces to the system
\begin{eqnarray}
\alpha,\beta_2&\not\equiv&0\pmod{p}\label{Vfirst}\\
a_1\alpha + a_3\gamma_1&\equiv&\beta_2\alpha c_1 + \beta_3\alpha d_1\pmod{p}
\nonumber\\
a_3 \gamma_2&\equiv&\beta_2\lambda_1\pmod{p}\label{Vthird}\\
a_3&\equiv&\beta_2\alpha c_3\pmod{p}\label{Vfourth}\\ 
b_1&\equiv&d_1\pmod{p}\nonumber
\end{eqnarray}
\noindent \eqref{Vthird} can be ignored since $\lambda_1$
can always be chosen to solve this congruence given solutions to the others.
\eqref{Vfirst} and \eqref{Vfourth} imply that $\sigma \sim \tau$ only if 
$a_3 \equiv c_3 \equiv 0 \pmod{p}$ or $a_3,c_3 \not\equiv 0 \pmod{p}$, and so 
this gives us the following two mutually
exclusive cases.
\begin{enumerate}\item $a_3 \equiv c_3 \equiv 0 \pmod{p}$\\
If $b_1 \equiv d_1 \equiv 0 \pmod{p}$ then the system 
can be solved if and only if $a_1 \equiv c_1 \equiv 0 \pmod{p}$ or
$a_1,c_1 \not\equiv 0 \pmod{p}$. If $b_1 \equiv d_1 \not\equiv 0 \pmod{p}$ then
the system can always be solved. We therefore have $p+1$ distinct $\sim$-classes
in this case with representatives
\begin{equation}
\label{case1Vreps}
\left(
\begin{array}{ccc}
0&0&0\\
p&0&0\\
0&0&0
\end{array}
\right)
\quad
\left(
\begin{array}{ccc}
0&0&0\\
0&0&0\\
p\epsilon&0&0
\end{array}
\right)
\quad\epsilon=0,1,\ldots,p-1.
\end{equation}
\item $a_3,c_3 \not\equiv 0 \pmod{p}$\\
In this case the system is always soluble and so we have $p$ distinct
$\sim$-classes with representatives
\begin{equation}
\label{case2Vreps}
\left(
\begin{array}{ccc}
0&0&0\\
0&0&0\\
p\epsilon&0&0
\end{array}
\right)
\quad\epsilon=0,1,\ldots,p-1.
\end{equation}
\end{enumerate}
\noindent Therefore for the Lie ring $V$ there are $2p + 1$ $\sim$-classes
with representatives given in \eqref{case1Vreps} and \eqref{case2Vreps}. 

\subsubsection{The Lie ring $W$.}

We proceed as for $V$. First observe that any nilpotent element of 
${\rm Der}(W)_z$ can be represented in the form \eqref{nilpotentform}, and that
any such matrix 
$\kappa$ satisifies $[u_1,u_2]\kappa = [u_1\kappa,u_2] + [u_1,u_2\kappa]$.
Therefore since the inner derivations of $W$ are the matrices
\begin{equation*}
\left(
\begin{array}{ccc}
0&0&0\\
p{\mathbb Z}/p^2{\mathbb Z}&0&0\\
p{\mathbb Z}/p^2{\mathbb Z}&0&0
\end{array}
\right),
\end{equation*}

\noindent any nilpotent element of ${\rm Der}(W)_z$ is $\sim$-equivalent to a 
derivation of the form
\begin{equation}
\label{Wders}
\left(
\begin{array}{ccc}
0&0&0\\
a_1&p a_2&a_3\\
0&p b_2&0
\end{array}
\right).
\end{equation}

\noindent Now let $\pi$ be a matrix of the form \eqref{autform} where the 
diagonal
entries are coprime to $p$. Then the condition $\pi\in{\rm Aut}_{\rm Lie}(W)$ is
given by $\alpha pz = \beta_2\gamma_3 p z$, and this is equivalent to
$\alpha \equiv \beta_2 \gamma_3 \pmod{p}$. We now determine the conditions
for $\sim$-equivalence between derivations of the form \eqref{Wders}.
So letting 
\begin{equation*}
\sigma = \left(
\begin{array}{ccc}
0&0&0\\
a_1&pa_2&a_3\\
0&pb_2&0
\end{array}
\right)
\quad{\rm and}\quad
\tau = \left(
\begin{array}{ccc}
0&0&0\\
c_1&pc_2&c_3\\
0&pd_2&0
\end{array}
\right),
\end{equation*}
\noindent we see that $\sigma \sim \tau$ if and only if there exist matrices
\[
m = \left(
\begin{array}{ccc}
0&0&0\\
p\lambda_1&0&0\\
p\lambda_2&0&0
\end{array}
\right)
\quad{\rm and}\quad
\pi = \left(
\begin{array}{ccc}
\alpha&0&0\\
\beta_1&\beta_2&\beta_3\\
p\gamma_1&p\gamma_2&\gamma_3
\end{array}
\right),
\]
with $\alpha,\beta_2,\gamma_3\not\equiv 0\pmod{p},
 \alpha\equiv\beta_2\gamma_3\pmod{p}$, and such that 
$\sigma \pi = \pi (m + \alpha \tau)$. This is equivalent to the following 
system of congruences having a solution for the entries of $\pi$ and $m$
\begin{eqnarray}
\alpha,\beta_2,\gamma_3&\not\equiv&0\pmod{p}\label{Wfirst}\\
\alpha&\equiv&\beta_2\gamma_3\pmod{p}\label{Wsecond}\\
a_1\alpha + p(a_2\beta_1 + a_3\gamma_1)&\equiv&\beta_2\alpha c_1 + 
p(\beta_2\lambda_1 + \beta_3\lambda_2)\pmod{p^2}\label{Wthird}\\
p(a_2\beta_2 + a_3\gamma_2)&\equiv&p\alpha(\beta_2 c_2 + \beta_3 d_2)
\pmod{p^2}\label{Wfourth}\\
a_3\gamma_3&\equiv&\beta_2\alpha c_3\pmod{p}\label{Wfifth}\\
pb_2\beta_1&\equiv&p(\gamma_2\alpha c_1 + \gamma_3\lambda_2)\pmod{p^2}
\label{Wsixth}\\
pb_2\beta_2&\equiv&p\alpha d_2\gamma_3\pmod{p^2}\label{Wseventh}
\end{eqnarray}
\noindent \eqref{Wfirst},\eqref{Wsecond} and \eqref{Wthird} show that
this system can be solved only if $a_1\equiv c_1\equiv 0\pmod{p}$ or 
$a_1,c_1\not\equiv\pmod{p}$. We therefore have two mutually exclusive cases.
\begin{enumerate}
\item $a_1\equiv c_1\equiv 0\pmod{p}$.\\
Choosing $a_1^\prime,c_1^\prime$ so that $pa_1^\prime\equiv a_1 \pmod{p^2}$
and $pc_1^\prime\equiv c_1\pmod{p}$, we see that in this case, 
$\sigma \sim \tau$ if and
only if the system consisting of \eqref{Wfirst},\eqref{Wsecond} together with 
\begin{eqnarray}
a_1^\prime\alpha + a_2\beta_1 + a_3\gamma_1&\equiv&\beta_2\alpha c_1^\prime +
\beta_2\lambda_1 + \beta_3\lambda_2\pmod{p}\nonumber\\
a_2\beta_2 + a_3\gamma_2&\equiv&\alpha\beta_2c_2 + \beta_3d_2\pmod{p}
\label{Wonefirst}\\
a_3\gamma_3&\equiv&\beta_2\alpha c_3\pmod{p}\label{Wonesecond}\\
b_2\beta_1&\equiv&\gamma_3\lambda_2\pmod{p}\nonumber\\
b_2\beta_2&\equiv&\alpha d_2\gamma_3\pmod{p}\label{Wonethird}
\end{eqnarray}
\noindent has a solution. The congruences involving $\lambda_1$ and $\lambda_2$
can be ignored since they can always be solved given solutions to the others,
and so the system reduces to \eqref{Wfirst},\eqref{Wsecond},\eqref{Wonefirst},
\eqref{Wonesecond} and \eqref{Wonethird}. We have the following four mutually
exclusive subcases.
\begin{enumerate}
\item $a_3\equiv c_3 \equiv 0 \pmod{p}$ and $b_2\equiv d_2\equiv 0 \pmod{p}$.\\
\eqref{Wonefirst},\eqref{Wonesecond} and \eqref{Wonethird} reduce to the single
congruence $a_2\beta_2\equiv\alpha\beta_2 c_2\pmod{p}$, and the system 
consisting of this congruence together with
\eqref{Wfirst} and \eqref{Wsecond} has a solution if and only if
$a_2\equiv c_2\equiv 0\pmod{p}$ or $a_2,c_2\not\equiv 0\pmod{p}$. We have two
$\sim$-classes with representatives
\[
\left(
\begin{array}{ccc}
0&0&0\\
0&0&0\\
0&0&0\\
\end{array}
\right)
\quad{\rm and}\quad
\left(
\begin{array}{ccc}
0&0&0\\
0&p&0\\
0&0&0\\
\end{array}
\right).
\]
\item $a_3\equiv c_3\equiv 0\pmod{p}$ and $b_2,d_2 \not\equiv 0\pmod{p}$.\\
From \eqref{Wfirst},\eqref{Wsecond} and \eqref{Wonethird}, we see that
a necessary condition for $\sigma \sim \tau$ in this case is for 
$b_2 \equiv d_2\gamma_3^2 \pmod{p}$ to have a solution for $\gamma_3$, i.e.\ 
$b_2$ and $d_2$
must have the same quadratic character modulo $p$. It is easy to see that this
condition is sufficient also, and so letting $\nu$ denote a non-quadratic
residue modulo $p$, we have two $\sim$-classes with representatives
\[
\left(
\begin{array}{ccc}
0&0&0\\
0&0&0\\
0&p&0\\
\end{array}
\right)
\quad{\rm and}\quad
\left(
\begin{array}{ccc}
0&0&0\\
0&0&0\\
0&p\nu&0\\
\end{array}
\right).
\]
\item $a_3,c_3\not\equiv 0\pmod{p}$ and $b_2\equiv d_2\equiv 0\pmod{p}$.\\
In this case, a necessary and sufficient condition is for 
$a_3\equiv c_3\beta_2^2 \pmod{p}$ to have a solution for $\beta_2$, and so 
letting $\nu$ denote a non-quadratic residue modulo $p$, we have two
$\sim$-classes with representatives
\[
\left(
\begin{array}{ccc}
0&0&0\\
0&0&1\\
0&0&0\\
\end{array}
\right)
\quad{\rm and}\quad
\left(
\begin{array}{ccc}
0&0&0\\
0&0&\nu\\
0&0&0\\
\end{array}
\right).
\]
\item $a_3,c_3\not\equiv 0\pmod{p}$ and $b_2,d_2\not\equiv 0\pmod{p}$.\\
In this case, necessary and sufficient conditions for $\sigma\sim\tau$ are that 
$a_3\equiv c_3\beta_2^2\pmod{p}$ and $b_2\equiv d_2\gamma_3^2\pmod{p}$ have
solutions. Therefore letting $\nu$ denote a non-quadratic residue modulo $p$,
we have four $\sim$-classes with representatives
\[
\qquad\qquad\left(
\begin{array}{ccc}
0&0&0\\
0&0&1\\
0&p&0\\
\end{array}
\right),
\left(
\begin{array}{ccc}
0&0&0\\
0&0&\nu\\
0&p&0\\
\end{array}
\right),
\left(
\begin{array}{ccc}
0&0&0\\
0&0&1\\
0&p\nu&0\\
\end{array}
\right),
\left(
\begin{array}{ccc}
0&0&0\\
0&0&\nu\\
0&p\nu&0\\
\end{array}
\right).
\]
\end{enumerate}
\noindent So summarising, in case 1 we have 10 distinct $\sim$-classes. We now
deal with case 2.
\item $a_1,c_1\not\equiv 0\pmod{p}$.\\
In this case, \eqref{Wfirst} and \eqref{Wthird} imply that we must have
$\beta_2\equiv a_1c_1^{-1}\pmod{p}$ in order for
$\sigma\sim\tau$. Now if we have a solution to the system consisting of 
\eqref{Wfirst},\eqref{Wsecond},\eqref{Wfourth},\eqref{Wfifth} and 
\eqref{Wseventh}
(i.e.\ the congruences not involving $\lambda_1,\lambda_2$), together with
this constraint on $\beta_2$, then $\lambda_2$ can be chosen to solve
\eqref{Wsixth} (since $\gcd(\gamma_3,p)=1$). Then choosing 
$\lambda_1\equiv 0\pmod{p}$ we can lift $\beta_2$ to give a solution to
\eqref{Wthird} (since $\gcd(\alpha c_1,p)=1$). Therefore, the system reduces
to \eqref{Wfirst},\eqref{Wsecond},\eqref{Wfifth}, together with
\begin{eqnarray}
\beta_2&\equiv&a_1c_1^{-1}\pmod{p}\label{Wsecondone}\\
a_2\beta_2 + a_3\gamma_2&\equiv&c_2\alpha\beta_2 + d_2\alpha\beta_3\pmod{p}
\label{Wsecondtwo}\\
b_2\beta_2&\equiv&d_2\alpha\gamma_3\pmod{p}\label{Wsecondthree}
\end{eqnarray} 
\noindent Now since we must have $\alpha,\beta_2,\gamma_3$ coprime to $p$,
\eqref{Wfifth} and \eqref{Wsecondthree} imply that we have the following
four mutually exclusive subcases.
\begin{enumerate}
\item $a_3\equiv c_3\equiv 0\pmod{p}$ and $b_2\equiv d_2\equiv 0\pmod{p}$.\\
\eqref{Wsecondone} and \eqref{Wsecondtwo} imply that $\sigma\sim\tau$ only if
$a_2\equiv c_2\equiv 0\pmod{p}$ or $a_2,c_2\not\equiv 0\pmod{p}$, and either of
these is seen to be sufficient as well, so that in this subcase we have two 
$\sim$-classes with representatives
\[
\left(
\begin{array}{ccc}
0&0&0\\
1&0&0\\
0&0&0\\
\end{array}
\right)
\quad{\rm and}\quad
\left(
\begin{array}{ccc}
0&0&0\\
1&p&0\\
0&0&0\\
\end{array}
\right).
\]
\item $a_3\equiv c_3\equiv 0\pmod{p}$ and $b_2,d_2\not\equiv 0\pmod{p}$.
Using \eqref{Wsecond} to substitute for $\alpha$ in \eqref{Wsecondthree}, we 
see that
a necessary condition for $\sigma\sim\tau$ is that $b_2$ and $d_2$ have the
same quadratic character modulo $p$. This condition is easily seen to be
sufficient also, and so denoting by $\nu$ a non-quadratic residue modulo $p$,
we have two $\sim$-classes with representatives
\[
\left(
\begin{array}{ccc}
0&0&0\\
1&0&0\\
0&p&0\\
\end{array}
\right)
\quad{\rm and}\quad
\left(
\begin{array}{ccc}
0&0&0\\
1&0&0\\
0&p\nu&0\\
\end{array}
\right).
\]
\item $a_3,c_3\not\equiv 0\pmod{p}$ and $b_2\equiv d_2\equiv 0\pmod{p}$.\\
\eqref{Wsecondtwo} together with \eqref{Wfirst} and \eqref{Wsecond} always
has a solution, and so eliminating $\alpha$ and $\beta_2$ from \eqref{Wfifth}
we see that $\sigma\sim\tau$ if and only if 
$a_3c_3^{-1}\equiv(a_1c_1^{-1})^2\pmod{p}$. Now for each pair $(a_1,a_3)$ in
$({\mathbb Z}/p{\mathbb Z})^\times \times ({\mathbb Z}/p{\mathbb Z})^\times$
there are $((p-1)/2)\times 2$ pairs $(c_1,c_3)$ in
$({\mathbb Z}/p{\mathbb Z})^\times \times ({\mathbb Z}/p{\mathbb Z})^\times$
satisfying this congruence. Therefore there are
$\vert ({\mathbb Z}/p{\mathbb Z})^\times \times 
({\mathbb Z}/p{\mathbb Z})^\times\vert / p-1$ distinct $\sim$-classes in this
subcase. Letting $h$ be a primitive element modulo $p$ we have the following 
representatives.
\[
\left(
\begin{array}{ccc}
0&0&0\\
h^r&0&1\\
0&0&0\\
\end{array}
\right)
\quad
\left(
\begin{array}{ccc}
0&0&0\\
h^r&0&h\\
0&0&0\\
\end{array}
\right)
\quad
{\rm where}\ r = 1,\ldots,(p-1)/2.
\]
\item $a_3,c_3\not\equiv 0\pmod{p}$ and $b_2,d_2\not\equiv 0\pmod{p}$.\\
Using \eqref{Wfirst},\eqref{Wsecond} and \eqref{Wsecondone} we can eliminate 
$\alpha$
and $\beta_2$ from \eqref{Wfifth} and \eqref{Wsecondthree} to obtain
\begin{eqnarray}
a_3c_3^{-1}&\equiv&(a_1c_1^{-1})^2 \pmod{p}\label{Wsecondeighth}\\
b_2&\equiv&d_2\gamma_3^2\pmod{p}\nonumber
\end{eqnarray}
A solution to these two congruences guarantees the existence of a solution to 
the whole
system, and so $\sigma\sim\tau$ in this subcase if and only if $b_2$ and $d_2$
have the same quadratic character modulo $p$, and \eqref{Wsecondeighth} holds.
By a similar discussion to the previous subcase, we see that there are
$2(p-1)$ $\sim$-classes. Letting $h$ be a primitive element modulo $p$ (so 
that, in particular, $h$ is a non-quadratic residue modulo $p$), we have the
following representatives where $r=1,\ldots,(p-1)/2$
\begin{equation*}
\qquad\qquad\left(
\begin{array}{ccc}
0&0&0\\
h^r&0&1\\
0&p&0
\end{array}
\right)
\left(
\begin{array}{ccc}
0&0&0\\
h^r&0&h\\
0&p&0
\end{array}
\right)
\left(
\begin{array}{ccc}
0&0&0\\
h^r&0&1\\
0&ph&0
\end{array}
\right)
\left(
\begin{array}{ccc}
0&0&0\\
h^r&0&h\\
0&ph&0
\end{array}
\right)
\end{equation*}
\end{enumerate}
\noindent So summarising, in case 1 we have shown that there are
$3(p-1)+4=3p+1$ distinct $\sim$-classes.
\end{enumerate}
\noindent Therefore for the Lie ring $W$, we have shown that there are
$3p+11$ distinct $\sim$-classes.

\subsubsection{The Lie ring $X$.}
Since $X$ is the Abelian Lie ring on $A$, all inner derivations are 0 and the 
nilpotent elements of ${\rm Der}(X)_z$ are precisely the matrices of the form
\eqref{nilpotentform}. Also, the automorphisms of $X$ fixing 
$\langle z \rangle$ set-wise are
the matrices represented as \eqref{autform} with the diagonal entries coprime
to $p$. So if we let $\sigma$ and $\tau$ be two elements of ${\rm Der}(X)_z$
where 
\[
\sigma=\left(
\begin{array}{ccc}
0&0&0\\
a_1&pa_2&a_3\\
pb_1&pb_2&0
\end{array} 
\right)
\quad{\rm and}\quad
\tau=\left(
\begin{array}{ccc}
0&0&0\\
c_1&pc_2&c_3\\
pd_1&pd_2&0
\end{array}
\right)
,
\]
then $\sigma\sim\tau$ if and only if there exists a matrix
\[
\pi=\left(
\begin{array}{ccc}
\alpha&0&0\\
\beta_1&\beta_2&\beta_3\\
p\gamma_1&p\gamma2&\gamma_3
\end{array}
\right)
\quad{\rm with}\ \gcd(\alpha,p)=\gcd(\beta_2,p)=1,
\]
and such that $\sigma\pi = \pi\alpha\tau$. This gives the following system
\begin{eqnarray}
\alpha,\beta_2,\gamma_3&\not\equiv&0\pmod{p}\label{coprime}\\
a_1\alpha + p(a_2\beta_1 + a_3\gamma_1)&\equiv&\beta_2\alpha c_1 + 
p\beta_3\alpha d_1 \pmod{p^2}\label{Xone}\\
p(a_2\beta_2 + a_3\gamma_2)&\equiv&p(\alpha\beta_2 c_2 + \alpha\beta_3 d_2)
\pmod{p^2}\label{Xtwo}\\
a_3\gamma_3&\equiv&\beta_2\alpha c_3\pmod{p}\label{Xthree}\\
p(b_1\alpha + b_2\beta_1)&\equiv&p(\alpha\gamma_2 c_1 + \alpha\gamma_3 d_1)
\pmod{p^2}\label{Xfour}\\
pb_2\beta_2&\equiv&pd_2\alpha\gamma_3\pmod{p^2}\label{Xfive}
\end{eqnarray}
\eqref{coprime} and \eqref{Xone} imply that we must have 
$a_1\equiv c_1\equiv 0\pmod{p}$ or $a_1,c_1\not\equiv 0\pmod{p}$ in order to
solve this system, and so we 
have two mutually exclusive cases.
\begin{enumerate}
\item $a_1\equiv c_1\equiv 0\pmod{p}$.\\
Writing $a_1\equiv pa_1^\prime \pmod{p}$ and $c_1\equiv pc_1^\prime \pmod{p}$,
we have $\sigma\sim\tau$ if and only if the following system is soluble
together with \eqref{coprime} and \eqref{Xthree}
\begin{eqnarray}
a_1^\prime\alpha + a_2\beta_1 + a_3\gamma_1&\equiv&c_1^\prime\alpha\beta_2 +
d_1\alpha\beta_3\pmod{p}\label{Xoneone}\\
a_2\beta_2 + a_3\gamma_2&\equiv&c_2\alpha\beta_2 + d_2\alpha\beta_3 \pmod{p}
\label{Xonetwo}\\
b_1\alpha + b_2\beta_1&\equiv&d_1\alpha\gamma_3\pmod{p}\label{Xonethree}\\
b_2\beta_2&\equiv&d_2\alpha\gamma_3\pmod{p}\label{Xonefour}
\end{eqnarray}
\eqref{coprime} together with \eqref{Xthree} and \eqref{Xonefour} give us four
mutually exclusive subcases.
\begin{enumerate}
\item $a_3\equiv c_3\equiv 0 \pmod{p}$ and $b_2\equiv d_2\equiv 0 \pmod{p}$.\\
The system we consider here is \eqref{coprime} together with
\begin{eqnarray}
a_1^\prime\alpha + a_2\beta_1&\equiv&c_1^\prime\alpha\beta_2 + d_1\alpha\beta_3
\pmod{p}\label{Xoneoneone}\\
a_2&\equiv&c_2\alpha\pmod{p}\label{Xoneonetwo}\\
b_1&\equiv&d_1\gamma_3\pmod{p}\label{Xoneonethree}
\end{eqnarray}
\eqref{coprime} together with \eqref{Xoneonetwo} and \eqref{Xoneonethree} give 
us four mutually exclusive subcases.
\begin{enumerate}
\item $a_2\equiv c_2\equiv 0\pmod{p}$ and $b_1\equiv d_1\equiv 0\pmod{p}$.\\
The condition for $\sigma\sim\tau$ is given by \eqref{coprime} and
$a_1^\prime\equiv c_1^\prime\beta_2\pmod{p}$, and so we have two $\sim$-classes
(corresponding to $a_1^\prime\equiv c_1^\prime\equiv 0\pmod{p}$ and 
$a_1^\prime,c_1^\prime\not\equiv 0\pmod{p}$) with representatives
\[
\left(
\begin{array}{ccc}
0&0&0\\
0&0&0\\
0&0&0\\
\end{array}
\right)
\quad{\rm and}\quad
\left(
\begin{array}{ccc}
0&0&0\\
p&0&0\\
0&0&0\\
\end{array}
\right)
.
\]
\item $a_2\equiv c_2\equiv 0\pmod{p}$ and $b_1,d_1\not\equiv 0\pmod{p}$.\\
We have one $\sim$-class with representative
\[
\left(
\begin{array}{ccc}
0&0&0\\
0&0&0\\
p&0&0\\
\end{array}
\right).
\]
\item $a_2,c_2\not\equiv 0\pmod{p}$ and $b_1\equiv d_1\equiv 0\pmod{p}$.\\
We have one $\sim$-class with representative
\[
\left(
\begin{array}{ccc}
0&0&0\\
0&p&0\\
0&0&0\\
\end{array}
\right).
\]
\item $a_2,c_2\not\equiv 0\pmod{p}$ and $b_1,d_1\not\equiv 0\pmod{p}$.
$\alpha$ and $\gamma_3$ are uniquely defined by \eqref{Xoneonetwo} and
\eqref{Xoneonethree}, and then $\beta_1,\beta_2,\beta_3$ can be chosen to
solve \eqref{Xoneoneone}. So we have one $\sim$-class with representative
\[
\left(
\begin{array}{ccc}
0&0&0\\
0&p&0\\
p&0&0\\
\end{array}
\right).
\]
\end{enumerate}
\item $a_3\equiv c_3\equiv 0\pmod{p}$\\
Consider the matrix 
$S=\left(\begin{smallmatrix}c_1^\prime&d_1\\c_2&d_2\end{smallmatrix}\right)$,
and suppose first that $\det(S)\equiv 0\pmod{p}$ and $\sigma$ is the derivation
in the present case with 
$a_1^\prime\equiv a_2\equiv b_1\equiv 0\pmod{p}$ and $b_2\equiv 1\pmod{p}$.
Since $d_2\not\equiv 0\pmod{p}$, \eqref{Xoneone} is a multiple of 
\eqref{Xonetwo} and so can be ignored. Then using the fact that $d_2$ is
coprime to $p$, it is easy to see that $\sigma\sim\tau$. Therefore when
$\det(S)\equiv 0\pmod{p}$, we have shown that $\tau$ is $\sim$-equivalent to
\[
\rho_1=\left(
\begin{array}{ccc}
0&0&0\\
0&0&0\\
0&p&0\\
\end{array}
\right).
\]
Now suppose ${\rm det}(S)\not\equiv 0\pmod{p}$ and $\sigma$ is the derivation in
the present case with $a_1^\prime\equiv 0\pmod{p}$ and 
$a_2\equiv b_1\equiv b_2\equiv 1\pmod{p}$. Then using \eqref{Xonethree} to
substitute
for $\beta_1$ in \eqref{Xoneone}, and using \eqref{Xonefour} to substitute for 
$\beta_2$ in
\eqref{Xonetwo}, we can replace \eqref{Xoneone} and \eqref{Xonetwo} with
\begin{eqnarray}
d_1\gamma_3 - 1&\equiv&c_1^\prime\beta_2 + d_1\beta_3\nonumber\\
d_2\gamma_3&\equiv&c_2\beta_2 + d_2\beta_3\nonumber
\end{eqnarray}
For any value of $\gamma_3$ these have a unique solution for 
$(\beta_2,\beta_3)$ (and $\beta_2\not\equiv 0\pmod{p}$ since 
$d_1\not\equiv 0\pmod{p}$). Therefore $\tau\sim\rho_2$ where
\[
\rho_2=\left(
\begin{array}{ccc}
0&0&0\\
0&p&0\\
p&p&0\\
\end{array}
\right).
\]
To see that $\rho_1\not\sim\rho_2$, observe that the derived subring of 
$X(4,w_4,\rho_1,z)$ has order $p$ whereas the derived subring of 
$X(4,w_4,\rho_2,z)$ has order $p^2$ (this notation was defined in 
\lemref{construct}). Therefore we have two $\sim$-classes with
representatives $\rho_1,\rho_2$. 
\item $a_3,c_3\not\equiv 0\pmod{p}$ and $b_2\equiv d_2\equiv 0\pmod{p}$.\\
Given solutions to \eqref{Xthree} and \eqref{Xonethree} satisfying
\eqref{coprime}, we can easily obtain solutions to \eqref{Xoneone},
\eqref{Xonetwo} and \eqref{Xonefour}. Now \eqref{Xonethree} is soluble if and
only if $b_1\equiv d_1\equiv 0\pmod{p}$ or $b_1,d_1\not\equiv 0\pmod{p}$, and 
so we have two $\sim$-classes with representatives
\[
\left(
\begin{array}{ccc}
0&0&0\\
0&0&1\\
0&0&0\\
\end{array}
\right)
\quad{\rm and}\quad
\left(
\begin{array}{ccc}
0&0&0\\
0&0&1\\
p&0&0\\
\end{array}
\right)
.
\]
\item $a_3,c_3\not\equiv 0\pmod{p}$ and $b_2,d_2\not\equiv 0\pmod{p}$.\\
Given solutions to \eqref{Xthree} and \eqref{Xonefour} satisfying 
\eqref{coprime}, we can easily find solutions to \eqref{Xoneone},
\eqref{Xonetwo} and \eqref{Xonethree}. Now using \eqref{Xonefour} to eliminate
$\alpha$ from \eqref{Xthree} we obtain 
$a_3c_3^{-1}d_2b_2^{-1}\equiv\beta_2^2\gamma_3^{-2}\pmod{p}$, and this is
soluble if and only if $a_3c_3$ and $d_2b_2$ have the same quadratic character
modulo $p$. Letting $\nu$ denote a non-quadratic residue modulo $p$, we have
two $\sim$-classes with representatives
\[
\left(
\begin{array}{ccc}
0&0&0\\
0&0&1\\
0&p&0\\
\end{array}
\right)
\quad{\rm and}\quad
\left(
\begin{array}{ccc}
0&0&0\\
0&0&1\\
0&p\nu&0\\
\end{array}
\right).
\]
\end{enumerate}
So summarising, in case 1 there are 11 distinct $\sim$-classes.
\item $a_1,c_1\not\equiv 0\pmod{p}$.\\
From \eqref{coprime} and \eqref{Xone} we see that a necessary condition for
$\sigma\sim\tau$ is $\beta_2\equiv a_1c_1^{-1}\pmod{p}$. With this constraint
on $\beta_2$, given any solution to \eqref{Xtwo}--\eqref{Xfive} we can lift
$\beta_2$ to solve \eqref{Xone}. Therefore the system we consider is
\eqref{coprime},\eqref{Xthree} and
\begin{eqnarray}
\beta_2&\equiv&a_1c_1^{-1}\pmod{p}\label{Xtwoone}\\
a_2\beta_2+a_3\gamma_2&\equiv&\alpha\beta_2c_2+\alpha\beta_3d_2\pmod{p}
\label{Xtwotwo}\\
b_1\alpha+b_2\beta_1&\equiv&\alpha\gamma_2c_1+\alpha\gamma_3d_1\pmod{p}
\label{Xtwothree}\\
b_2\beta_2&\equiv&d_2\alpha\gamma_3\pmod{p}\label{Xtwofour}
\end{eqnarray}
From \eqref{Xthree} and \eqref{Xtwofour} we have four mutually exclusive
subcases.
\begin{enumerate}
\item $a_3\equiv c_3\equiv 0\pmod{p}$ and $b_2\equiv d_2\equiv 0\pmod{p}$.\\
\eqref{Xtwotwo} reduces to $a_2\equiv c_2\alpha\pmod{p}$ and so 
$\sigma\sim\tau$ only if $a_2\equiv c_2\equiv 0\pmod{p}$ or 
$a_2,c_2\not\equiv 0\pmod{p}$. Either of these is sufficient also, and so we
have two $\sim$-classes with representatives
\[
\left(
\begin{array}{ccc}
0&0&0\\
1&0&0\\
0&0&0\\
\end{array}
\right)
\quad{\rm and}\quad
\left(
\begin{array}{ccc}
0&0&0\\
1&p&0\\
0&0&0\\
\end{array}
\right).
\]
\item $a_3\equiv c_3\equiv 0\pmod{p}$ and $b_2,d_2\not\equiv 0\pmod{p}$.\\
Choosing $\gamma_3\not\equiv 0\pmod{p}$ arbitrarily, \eqref{Xtwoone} and
\eqref{Xtwofour} define $\alpha$, and then $\gamma_2,\beta_3$ can be chosen
to give a solution to the other congruences. Therefore there is one 
$\sim$-class with representative
\[
\left(
\begin{array}{ccc}
0&0&0\\
1&0&0\\
0&p&0\\
\end{array}
\right).
\]
\item $a_3,c_3\not\equiv 0\pmod{p}$ and $b_2\equiv d_2\equiv 0\pmod{p}$.\\
First consider the situation where $a_2\equiv c_2\equiv 0\pmod{p}$. Then from
\eqref{Xtwotwo} we see that $\gamma_2\equiv 0\pmod{p}$, and so from
\eqref{Xtwothree} we see that $\sigma\sim\tau$ only if 
$b_1\equiv d_1\equiv 0\pmod{p}$ or $b_1,d_1\not\equiv 0\pmod{p}$. Therefore
there are at least two $\sim$-classes in this subcases with representatives
\[
\kappa_1=\left(
\begin{array}{ccc}
0&0&0\\
1&0&1\\
0&0&0\\
\end{array}
\right)
\quad{\rm and}\quad
\kappa_2=\left(
\begin{array}{ccc}
0&0&0\\
1&p&0\\
p&0&0\\
\end{array}
\right).
\]
Now let $\tau$ be an arbitrary derivation in this subcase and suppose first
that $c_3d_1+c_2c_1\equiv 0\pmod{p}$. Then letting $\sigma=\kappa_1$ (i.e.\
$a_1\equiv a_3\equiv 1\pmod{p}$ and $a_2\equiv b_1\equiv 0\pmod{p}$) it is
straightforward to see that $\tau\sim\sigma$. Now suppose that
$c_3d_1+c_2c_1\not\equiv 0\pmod{p}$ and let $\sigma=\kappa_2$ (i.e.\
$a_1\equiv a_3\equiv b_1\equiv 1\pmod{p}$ and $a_2\equiv 0\pmod{p}$). The
condition for $\tau\sim\sigma$ is given by \eqref{coprime},\eqref{Xthree},
\eqref{Xtwoone} and 
\begin{eqnarray}
\gamma_2&\equiv&\alpha\beta_2c_2\pmod{p}\label{Xtwothreegamma2}\\
c_1\gamma_2 + d_1\gamma_3&\equiv&1\pmod{p}\label{Xtwothreeone}
\end{eqnarray}
Using \eqref{Xthree} and \eqref{Xtwothreegamma2} to substitute for $\gamma_3$ 
and $\gamma_2$ (respectively) in \eqref{Xtwothreeone}, 
we see that \eqref{Xtwothreeone} is equivalent to 
$\alpha^{-1}\equiv\beta_2(c_3d_1+c_2c_1)\pmod{p}$. Since \eqref{Xtwoone}
ensures that the right-hand side of this is coprime to $p$, it follows that
the system always has a solution. Therefore we have exactly two $\sim$-classes 
with in this subcase with representatives $\kappa_1,\kappa_2$.
\item $a_3,c_3\not\equiv 0\pmod{p}$ and $b_2,d_2\not\equiv 0\pmod{p}$.\\
\eqref{Xtwotwo} and \eqref{Xtwothree} can be ignored ($\beta_1$ and
$\beta_3$ can be chosen to solve them given solutions to the others), and
using \eqref{Xtwofour} to substitute for $\alpha$ in \eqref{Xthree}, we can 
replace \eqref{Xthree} with 
$a_3b_2^{-1}c_3^{-1}d_2\equiv\beta_2^2\gamma_3^{-2}\pmod{p}$. Therefore,
$\sigma\sim\tau$ only if $a_3b_2$ and $c_3d_2$ have the same quadratic
character modulo $p$, and this is seen to be sufficient also. So letting $\nu$
be a non-quadratic residue modulo $p$, we have two $\sim$-classes with 
representatives
\[
\left(
\begin{array}{ccc}
0&0&0\\
1&0&1\\
0&p&0\\
\end{array}
\right)
\quad{\rm and}\quad
\left(
\begin{array}{ccc}
0&0&0\\
1&0&1\\
0&p\nu&0\\
\end{array}
\right).
\]
\end{enumerate} 
So summarising, in case 2 we have shown that there are 7 $\sim$-classes.
\end{enumerate}
Therefore, for the Lie ring $X$ we have shown that there are $7+11=18$ 
$\sim$-classes.

\subsection{Calculating $\Psi_{n-3}^p(3)$ for $p\geq5$ and $n\geq 7$.}\mbox{}

From \propref{translate} and assuming the bound referred to in 
\remref{lieringbound}, we see that for $p\geq5$ and $n\geq7$, 
$\Psi_{n-3}^p(3)$ equals the number of isomorphism classes of nilpotent Lie 
ring structures on the Abelian group 
${\mathbb Z}/p^{n-3}{\mathbb Z} \oplus {\mathbb Z}/p^3{\mathbb Z}$. These
isomorphism classes are straightforward to calculate directly by using the
order of the derived subring as the principal invariant. We omit the complete
calculations but illustrate the general approach by showing that for a fixed
prime $p$ and integer $n\geq9$, there are three isomorphism classes where the
derived subring has order $p^3$.

So let $u,v$ be generators of ${\mathbb Z}/p^{n-3}{\mathbb Z}$ and
${\mathbb Z}/p^3{\mathbb Z}$ respectively, and suppose that $L$ is a nilpotent
Lie ring structure on 
$A = {\mathbb Z}/p^{n-3}{\mathbb Z} \oplus {\mathbb Z}/p^3{\mathbb Z}$
where $[u,v]$ has order $p^3$. Then 
$[u,v] = \alpha p^{n-6}u + \beta v$, and observe
that since $L$ is nilpotent we must have $\beta\equiv 0\pmod{p}$.
Hence $\alpha\not\equiv 0\pmod{p}$ (in order that $[L,L] = p^3$). 

\begin{enumerate}
\item $\beta \not\equiv 0\pmod{p^2}$. Write $\beta=p\gamma$ so that
$\gamma\not\equiv0\pmod{p}$, and let
$u^\prime=\gamma^{-1}u$ and $v^\prime=\alpha^{-1}v$ where 
$\gamma^{-1}\gamma\equiv1\pmod{p^3}$ and 
$\alpha^{-1}\alpha\equiv1\pmod{p^{n-3}}$. Then 
$(u^\prime,v^\prime)$ is a basis of $A$ and 
$[u^\prime,v^\prime] = p^{n-6}u^\prime + pv^\prime$. So $L$
is isomorphic to the Lie ring on $A$ where 
$[u,v] = p^{n-6}u + pv$. 
\item $\beta\equiv 0\pmod{p^2}$ and $\beta\not\equiv 0\pmod{p^3}$. Write
$\beta=p^2\gamma$ so that $\gamma\not\equiv0\pmod{p}$, and then letting
$u^\prime=\gamma^{-1}u$ and $v^\prime=\alpha^{-1}v$ shows that $L$ 
is isomorphic to the Lie ring on $A$ 
where $[u,v]=p^{n-6}u + p^2 v$.
\item $\beta\equiv 0\pmod{p^3}$. Letting 
$u^\prime=u$ and $v^\prime=\alpha^{-1} v^\prime$ shows that $L$ is 
isomorphic to the Lie ring on $A$ with $[u,v] = p^{n-6}u$.
\end{enumerate} 

\noindent Thus there are at most 3 isomorphism classes, but in case $i=1,2,3$, 
observe that the derived subring is contained in
$\agemo_i(A)$ but not in $\agemo_{i+1}(A)$ and so each case 
corresponds to a distinct isomorphism class. 

The calculations corresponding to the derived subring having order dividing
$p^2$ are similar and will be omitted. The treatment for $n=7$ and $n=8$ is
analogous. We now summarise the isomorphism classes which occur for any $p$ and
$n\geq7$ by giving the value of the Lie bracket $[u,v]$ (obviously, $u$ depends
on $p$ and $n$, and $v$ depends on $p$).

\begin{enumerate}
\item Derived subring of order $p^3$.

\begin{tabular}{lllll}
$n=7$&\quad&$pu$&&\\
$n=8$&\quad&$p^2u$&$p^2u+pv$\\
$n\geq9$&\quad&$p^{n-6}u+pv$&$p^{n-6}u+p^2v$&$p^{n-6}u$\\
\end{tabular}
\medskip
\item Derived subring of order $p^2$.

\begin{tabular}{lllll}
$n=7$&\quad&$p^2u$&$pv$&\\
$n\geq8$&\quad&$p^{n-5}u$&$p^{n-5}u+p^2v$&$pv$\\
\end{tabular}
\medskip
\item Derived subring of order 1 or $p$.

\begin{tabular}{lllll}
$n\geq7$&\quad&$p^{n-4}u$&$p^2v$&0\\
\end{tabular}

\end{enumerate}

\noindent From this list we obtain the expression for $\Psi_{n-3}^p(3)$ given
in \propref{decomposition}.
 
\bibliography{reg}
\bibliographystyle{plain}

\end{document}